%% file: Main.tex
\documentstyle[twoside,amssymb,12pt,thmsb,psfig,fullpage,sw20bams]{article}
\input{tcilatex}
\begin{document}

\title{Infinitesimal Thurston Rigidity and the Fatou-Shishikura Inequality}
\author{Adam Lawrence Epstein}
\date{\relax}
\maketitle

\thispagestyle{empty} \input{imsmark} \SBIMSMark{1999/1}{February 1999}{}

\section{\label{intro}Introduction}

Any rational map $f:{\Bbb P}^{1}\rightarrow {\Bbb P}^{1}$ of degree $D>1$
has infinitely many periodic points, by the Fundamental Theorem of Algebra.
The finiteness of the set of {\em nonrepelling} points is a cornerstone of
complex analytic dynamics. Recall that the cycle

\[
\langle x\rangle =\{x,...,f^{\kappa -1}(x)\}\text{ \quad is \quad }\left\{ 
\begin{array}{ll}
\text{{\em superattracting}} & \text{if \thinspace }\rho =0 \\ 
\text{{\em attracting}} & \text{if }0<|\rho |<1 \\ 
\text{{\em indifferent}} & \text{if }|\rho |=1 \\ 
\text{{\em repelling}} & \text{if }|\rho |>1,
\end{array}
\right. 
\]
where $\rho =(f^{\kappa })^{\prime }(x)$ is the corresponding {\em %
eigenvalue. }An indifferent cycle is {\em rationally indifferent} if $\rho $
is a root of unity, and {\em irrationally indifferent }otherwise. The
assumption $D>1$ guarantees that every rationally indifferent cycle is {\em %
parabolic}, the first return to $x$ being expressible as 
\begin{equation}
\zeta \mapsto \rho (\zeta +\zeta ^{N+1}+\alpha \zeta ^{2N+1})+O(\zeta
^{2N+2})  \label{parnf}
\end{equation}
in a suitable local coordinate: if $\rho $ is a primitive $n$-root of unity
then $N=\nu n$ for some positive integer $\nu $ (see \cite{BH} or \cite{Mil}%
). The {\em Fatou-Shishikura Inequality} asserts that there are at most $%
2D-2 $ nonrepelling cycles, each parabolic cycle of the form (\ref{parnf})
counting as $\nu \geq 1$.

Here we present a new and independent proof of a refined Fatou-Shishikura
Inequality. The refinement concerns a more generous convention for counting
parabolic cycles: in terms of the normal form (\ref{parnf}) we associate to
each cycle $\langle x\rangle $ the quantity 
\[
\gamma _{\langle x\rangle }=\left\{ 
\begin{array}{ll}
0 & \text{ \thinspace if }\langle x\rangle \text{ is repelling or
superattracting} \\ 
1 & \text{ \thinspace if }\langle x\rangle \text{ is attracting or
irrationally indifferent} \\ 
\nu & \text{ \thinspace if }\langle x\rangle \text{ is {\em %
parabolic-repelling}\quad \thinspace }(\Re \beta >0) \\ 
\nu +1 & 
\begin{array}{l}
\text{if }\langle x\rangle \text{ is {\em parabolic-attracting } }(\Re \beta
<0) \\ 
\quad \quad \text{ }\,\text{or {\em parabolic-indifferent }}(\Re \beta =0),
\end{array}
\end{array}
\right. \text{ } 
\]
where $\beta =\frac{N+1}{2}-\alpha $ (in view of (\ref{iterate}) below, this
invariant behaves iteratively like $\frac{1}{\log \rho }$). Our count of
nonrepelling cycles of $f$ is $\gamma (f)=\sum_{\langle x\rangle \subset 
{\Bbb P}^{1}}\gamma _{\langle x\rangle }$ which {\em a priori} might be
infinite. We denote by $\delta (f)$ the number of infinite tails of critical
orbits; this quantity is certainly no greater than $2D-2$ (the number of
critical points), but if there are any critical orbit relations then it will
be smaller. Our refinement of the Fatou-Shishikura Inequality is the
following:

\begin{theorem}
\label{fs}Let $f:{\Bbb P}^{1}\rightarrow {\Bbb P}^{1}$ be a rational map of
degree $D>1$. Then $\gamma (f)\nolinebreak \leq \nolinebreak \delta (f).$
\end{theorem}

\smallskip

This formulation of the Fatou-Shishikura Inequality has the advantage that
the degree no longer explicitly appears, so that the appropriate extension
to transcendental maps is an assertion with content; the expanded account 
\cite{algdyn} of our argument contains a uniform treatment for all {\em %
finite type} analytic maps. We recover the usual formulation by observing
that there are at most $2D-2-\delta (f)$ superattracting cycles:

\begin{corollary}
\label{2d2}The number of superattracting, attracting, or indifferent cycles
is at most $2D-2$, for any rational map of degree $D>1.$
\end{corollary}

\smallskip

This result has a long history. Fatou \cite{fatou} and Julia \cite{Jul} both
proved that any attracting cycle must attract a critical point; in fact, a
parabolic cycle must have a critical point in each of the $\nu $ cycles of
petals (see \cite{BH} or \cite{Mil} for details). The relation between
critical points and irrationally indifferent cycles is rather more subtle.
It is fairly easy to show that any Cremer cycle lies in the postcritical
accumulation, and the same is true for Siegel disk boundaries; however, the
same critical point might well have an orbit which accumulates on several
such features, so this consideration does not even show that the set of
indifferent cycles is finite (Kiwi \cite{Kiwi} has recently circumvented
this difficulty in the polynomial case). Fatou instead found a perturbative
proof that the number of nonrepelling cycles is at most $4D-4$: loosely
speaking, half the indifferent cycles become attracting after a random
perturbation (see \cite{Mil}). Douady and Hubbard \cite{DH1} proved the
sharp bound for polynomials ($D-1$ nonrepelling cycles in ${\Bbb C}$) using
the theory of polynomial-like mappings, and finally Shishikura \cite{S}
proved the sharp bound $2D-2$ in complete generality. Shishikura employs
quasiconformal surgery to construct perturbations where all the nonrepelling
cycles become attracting. His discussion of the irrationally indifferent
case is especially delicate; in particular, the treatment of Cremer cycles
requires a careful comparison of asymptotics. Shishikura shows further,
using another beautiful surgery, that the same bound applies to the total
count of nonrepelling cycles augmented by twice the number of Herman ring
cycles.

Our proof of the Fatou-Shishikura Inequality is nonperturbative, and rather
more algebraic; in particular, we completely sidestep the classical
finiteness theorem for attracting cycles. The underlying mechanism is a
suitable extension of {\em Infinitesimal Thurston Rigidity}: the injectivity
of the linear operator $\nabla _{f}\nolinebreak =\nolinebreak I\nolinebreak
-\nolinebreak f_{*}$ on spaces of meromorphic quadratic differentials. More
precisely, the infinitesimal content of Thurston's Uniqueness Theorem (see 
\cite{DH2} and also \cite{Curt} for Global Rigidity and Thurston's much
harder Existence Theorem) is the assertion $f_{*}q\neq q$ for nonzero
quadratic differentials $q$ having at worst simple poles. The novel feature
of our extension is an allowance for multiple poles: for $q\in \ker \nabla
_{f}$ these are necessarily situated along cycles, so that the associated
infinities dynamically cancel, leaving finite residues whose signs reflect
the cycles' dynamical character (compare the proof \cite{J} of the Jenkins
General Coefficient Theorem). We address the issue of Herman rings in \cite
{algdyn}, and a planned sequel will discuss pertubative implications: we
will extend the considerations of \cite{thesis} to prove smoothness and
transversality for dynamically defined loci in parameter spaces.

\subsection*{Acknowledgments}

This work certainly never would have appeared without the moral support of
innumerable friends and colleagues. We owe particular and profound thanks to
John Hamal Hubbard for his unflagging patience and expository guidance. This
research was partially supported by NSF Grant DMS 9803242.

\section{\label{inj}Contraction, Injectivity and Finiteness}

Let ${\cal M}({\Bbb P}^{1})$ be the ${\Bbb C}$-linear space of all
meromorphic quadratic differentials $q$ on the Riemann sphere. We denote by $%
{\cal Q}({\Bbb P}^{1})$ the subspace consisting of all $q\in {\cal M}({\Bbb P%
}^{1})$ with at worst simple poles. Recall that $q\in {\cal Q}({\Bbb P}^{1})$
if and only if $\,||q||<\infty $, where $||q||\nolinebreak =\nolinebreak
\int_{{\Bbb P}^{1}}|q|$ is the total mass of the associated area form $|q|$;
see \cite{G} and also \cite{J} for the standard details. The quotient ${\cal %
M}({\Bbb P}^{1})/{\cal Q}({\Bbb P}^{1})$ is canonically isomorphic to ${\cal %
D}({\Bbb P}^{1})\nolinebreak =\nolinebreak \bigoplus_{x\in {\Bbb P}^{1}}%
{\cal D}_{x}({\Bbb P}^{1})$, where ${\cal D}_{x}({\Bbb P}^{1})$ is the space
of all {\em algebraic divergences} at $x$: polar parts, of order at most $-2$%
, of germs of meromorphic quadratic differentials. The algebraic divergence
of $q$ at $x$ is the corresponding class $[q]_{x}\nolinebreak \in
\nolinebreak {\cal D}_{x}({\Bbb P}^{1})$, and the total algebraic divergence
is $[q]=\sum_{x\in {\Bbb P}^{1}}[q]_{x}$. We write ${\cal M}({\Bbb P}^{1},A)$
for the subspace consisting of all $q\in {\cal M}({\Bbb P}^{1})$ whose poles
lie in a given set $A\subseteq {\Bbb P}^{1}$, and we denote by ${\cal Q}(%
{\Bbb P}^{1},A)$ the corresponding subspace of ${\cal Q}({\Bbb P}^{1})$. For
well-known cohomological reasons (see \nolinebreak \cite{nar}), $\dim {\cal Q%
}({\Bbb P}^{1},A)\nolinebreak =\nolinebreak \#A\nolinebreak -\nolinebreak 3$
so long as $\#A\nolinebreak \geq \nolinebreak 3$, and then ${\cal M}({\Bbb P}%
^{1},A)/{\cal Q}({\Bbb P}^{1},A)$ is canonically isomorphic to ${\cal D}(%
{\Bbb P}^{1},A)\nolinebreak =\nolinebreak \bigoplus_{x\in A}{\cal D}_{x}(%
{\Bbb P}^{1})$.

Recall that we may {\em pullback} any quadratic differential $q$ on ${\Bbb P}%
^{1}$ by any analytic map $\phi :U\rightarrow {\Bbb P}^{1}$ to obtain a
quadratic differential $\phi ^{*}q$ on $U$; the pullback of the associated
area form is $\phi ^{*}|q|=$ $|\phi ^{*}q|$. If $q$ is meromorphic at $\phi
(x)$ then $\phi ^{*}q$ is meromorphic at $x$: indeed, 
\begin{equation}
{\rm ord}_{x}\phi ^{*}q=\deg _{x}\phi \cdot ({\rm ord}_{\phi (x)}q+2)-2
\label{qorder}
\end{equation}
so a rational map $f:{\Bbb P}^{1}\rightarrow {\Bbb P}^{1}$ induces a
pullback operator $f^{*}\nolinebreak :\nolinebreak {\cal M}({\Bbb P}%
^{1})\nolinebreak \rightarrow \nolinebreak {\cal M}({\Bbb P}^{1})$ which
restricts to an endomorphism of ${\cal Q}({\Bbb P}^{1})$. Now consider the
corresponding {\em pushforward} operator $f_{*}:{\cal M}({\Bbb P}%
^{1})\rightarrow {\cal M}({\Bbb P}^{1})$: by definition, 
\[
f_{*}q=\sum_{\phi }\phi ^{*}q
\]
where $\phi $ ranges over the inverse branches of $f$, so that 
\begin{equation}
f_{*}f^{*}q=Dq  \label{pushpull}
\end{equation}
for a rational map of degree $D.$ As $||f_{*}q||=\int_{{\Bbb P}%
^{1}}|f_{*}q|\leq \int_{{\Bbb P}^{1}}f_{*}|q|=\int_{{\Bbb P}^{1}}|q|=||q||$
by the Triangle Inequality, the operator $f_{*}$ restricts to an
endomorphism of ${\cal Q}({\Bbb P}^{1})$; indeed, it follows directly from (%
\ref{qorder}) that 
\begin{equation}
{\rm ord}_{x}f_{*}q\geq \max_{w\in f^{-1}(x)}\left( \frac{{\rm ord}_{w}q+2}{%
\deg _{w}f}-2\right)   \label{ordpush}
\end{equation}
for each $x\in {\Bbb P}^{1}$, so ${\rm ord}_{x}f_{*}q\geq -1$ if ${\rm ord}%
_{w}q\geq -1$ for every $w\in f^{-1}(x)$. Furthermore, 
\mbox{${\rm ord}_{x}f_{*}q\geq 0$}
if ${\rm ord}_{w}q\geq 0$ for every such $w$, except
possibly when some $w$ is a critical point; we write $S(f)$ for the set of 
{\em critical values}, so that $P(f\nolinebreak )=\nolinebreak
\bigcup_{k=0}^{\infty }f^{k}(S(f))$ is the {\em postcritical set.} It
follows that 
\mbox{$f_{*}{\cal M}({\Bbb P}^{1},A) \subseteq 
{\cal M}({\Bbb P}^{1},f(A) \cup S(f))$}, and in particular 
\[
f_{*}{\cal Q}({\Bbb P}^{1},A)\subseteq {\cal Q}({\Bbb P}^{1},f(A)\cup S(f)),
\qquad \mbox{ for any\/} A\subseteq {\Bbb P}^{1};
\] similarly, 
\mbox{$f_{*}{\cal D}({\Bbb P}^{1},A) \subseteq {\cal D}({\Bbb P}^{1},f(A))$}
where 
\mbox{$f_{*}\negthinspace:{\cal D}({\Bbb P}^{1})\rightarrow 
{\cal D}({\Bbb P}^{1})$}
is the induced endomorphism.

Consider the ${\Bbb C}$-linear endomorphism $\nabla _{f}=I-f_{*}$ of ${\cal M%
}({\Bbb P}^{1})$. This operator restricts to an endomorphism of ${\cal Q}(%
{\Bbb P}^{1})$, so there is also an induced endomorphism $\nabla _{f}$ of $%
{\cal D}({\Bbb P}^{1})$. It follows from (\ref{ordpush}) that the space $%
{\cal D}(f)=\ker \nabla _{f}|_{{\cal D}({\Bbb P}^{1})}$ of {\em invariant
divergences} is computed cycle by cycle: to be precise, ${\cal D}%
(f)\nolinebreak =\nolinebreak \bigoplus_{\langle x\rangle \subset {\Bbb P}%
^{1}}{\cal D}_{\langle x\rangle }(f)$ where ${\cal D}_{\langle x\rangle
}(f)\nolinebreak =\nolinebreak \ker \nabla _{f}|_{{\cal D}({\Bbb P}%
^{1},\langle x\rangle )}$. Moreover, it suffices to compute these spaces for
fixed points: indeed, if $x$ is a point of period $\kappa $ then the
projection ${\cal D}({\Bbb P}^{1},\langle x\rangle )\nolinebreak \rightarrow
\nolinebreak {\cal D}_{x}({\Bbb P}^{1})$ restricts to an isomorphism ${\cal D%
}_{\langle x\rangle }(f)\nolinebreak \rightarrow \nolinebreak {\cal D}%
_{x}(f^{\kappa })$. We carry out this computation in Section \ref{local}: we
describe ${\cal D}_{\langle x\rangle }(f)$ in terms of the formal invariants 
$\rho $, $\nu $ and $\alpha $, then calculate the {\em dynamical residue } 
\[
{\cal D}_{\langle x\rangle }(f)\ni [q]_{\langle x\rangle }\longmapsto {\rm %
Res}_{\langle x\rangle }(f:q)\in {\Bbb R} 
\]
which measures the local creation or destruction of mass by $f_{*}$. The
relevant conclusions are summarized in:

\begin{proposition}
\label{dimd}The set ${\cal D}_{\langle x\rangle }^{\flat }(f)$ of all $%
[q]_{\langle x\rangle }\in {\cal D}_{\langle x\rangle }(f)$ with ${\rm Res}%
_{\langle x\rangle }(f:q)\leq 0$ is a ${\Bbb C}$-linear subspace of
dimension $\gamma _{\langle x\rangle }.$
\end{proposition}

\smallskip

We write ${\cal D}^{\flat }(f)$ for the subspace $\bigoplus_{\langle
x\rangle \subset {\Bbb P}^{1}}{\cal D}_{\langle x\rangle }^{\flat }(f)$ of $%
{\cal D}(f)$; for $A\subseteq {\Bbb P}^{1}$ we denote by ${\cal D}^{\flat
}(f,A)$ the corresponding subspace of ${\cal D}(f,A)=\bigoplus_{\langle
x\rangle \subseteq A}{\cal D}_{\langle x\rangle }(f)$. Consider the subspace 
\[
{\cal Q}^{\flat }(f)=\{q\in {\cal M}({\Bbb P}^{1}):[q]\in {\cal D}^{\flat
}(f)\} 
\]
of the ${\Bbb C}$-linear space ${\cal Q}(f)=\{q\in {\cal M}({\Bbb P}%
^{1}):[q]\in {\cal D}(f)\}=\nabla _{f}{}^{-1}{\cal Q}({\Bbb P}^{1})\supseteq
\ker \nabla _{f}.$ Our main result is the following:

\begin{proposition}
\label{main}Let $f:{\Bbb P}^{1}\rightarrow {\Bbb P}^{1}$ be a rational map
of degree $D>1$, and assume that $f$ is not a Latt\`{e}s example. Then $%
\nabla _{f}:{\cal Q}^{\flat }(f)\rightarrow {\cal Q}({\Bbb P}^{1})$ is
injective.
\end{proposition}

%\FRAME{ftbpFU}{3.4307in}{2.2762in}{0pt}{\Qcb{Net contraction}}{\Qlb{caricfig}%
%}{caricfig}{\special{language "Scientific Word";type
%"GRAPHIC";maintain-aspect-ratio TRUE;display "FULL";valid_file "F";width
%3.4307in;height 2.2762in;depth 0pt;original-width 0pt;original-height
%0pt;cropleft "0";croptop "1";cropright "1";cropbottom "0";filename
%'C:/papers/algdyn/caricfig.eps';file-properties "NPEU";}}

\begin{figure}[htbp]
\centerline{\psfig{figure=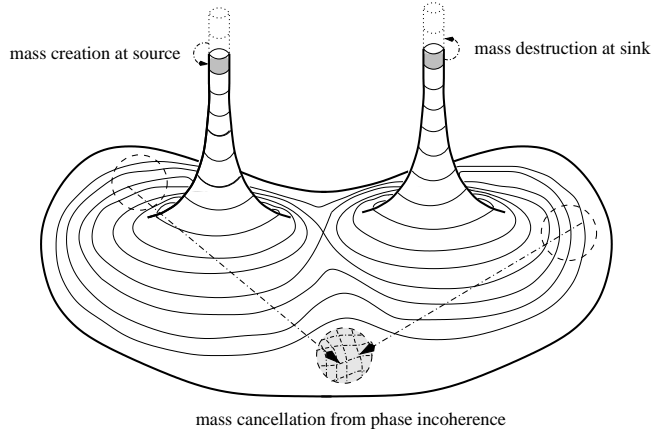,height=2.25in}}
\caption{Net contraction}
\label{caricfig}
\end{figure}

\smallskip We prove this algebraic assertion by combining two
measure-theoretic observations. The first generalizes the {\em Contraction
Principle} behind Thurston Rigidity: that cancellation due to phase
incoherence results in a well-defined decrease in mass. For nonintegrable
quadratic differentials this loss may be offset by the creation of new mass
at (parabolic\nolinebreak -\nolinebreak )repelling cycles or compounded by
the further destruction of mass at (parabolic\nolinebreak -\nolinebreak
)attracting cycles, and our second consideration is a {\em Balance Principle 
}which accordingly constains $\ker \nabla _{f}$. Figure \ref{caricfig} is a
caricature of this argument, which we present in Section \ref{global}; it is
not much harder to show that no eigenvalue of $f_{*}:{\cal M}({\Bbb P}%
^{1})\rightarrow {\cal M}({\Bbb P}^{1})$ lies on the unit circle \cite
{algdyn}.

\bigskip

\noindent {\bf Proof of Theorem \ref{fs}: }We may assume without loss of
generality that $f$ is not a Latt\`{e}s example, as such a map has only
repelling periodic points. Let $A$ be any finite set of the form $B\cup C$,
where $B\nolinebreak \subseteq \nolinebreak P(f)$ is an initial segment
including all critical orbit relations and $C$ consists of nonrepelling
cycles, and set $A^{+}=A\cup f(A)$. If $C\neq \varnothing $ then $\#A\geq 3$%
, so that 
\[
\begin{array}{ccccccccc}
0 & \longrightarrow & {\cal Q}({\Bbb P}^{1},A) & \longrightarrow & {\cal M}(%
{\Bbb P}^{1},A) & \longrightarrow & {\cal D}({\Bbb P}^{1},A) & 
\longrightarrow & 0 \\ 
&  & \text{ \quad }\downarrow {\small \nabla }_{f} &  & \text{ \quad }%
\downarrow {\small \nabla }_{f} &  & \text{ \quad }\downarrow {\small \nabla 
}_{f} &  &  \\ 
0 & \longrightarrow & {\cal Q}({\Bbb P}^{1},A^{+}) & \longrightarrow & {\cal %
M}({\Bbb P}^{1},A^{+}) & \longrightarrow & {\cal D}({\Bbb P}^{1},A^{+}) & 
\longrightarrow & 0
\end{array}
\]
is a commutative diagram of linear maps. As both rows are exact we may apply
the Serpent Lemma (an elementary diagram chase from homological algebra -
see \cite{nar}) to obtain an exact sequence 
\[
\ker \nabla _{f}\,|_{{\cal M}({\Bbb P}^{1},A)}\longrightarrow {\cal D}(f,A)%
\stackrel{\chi }{\longrightarrow }{\cal Q}({\Bbb P}^{1},A^{+})/\nabla _{f}%
{\cal Q}({\Bbb P}^{1},A)\text{.} 
\]
In view of Proposition \ref{main}, the restrictions $\chi |_{{\cal D}^{\flat
}(f,A)}$ and $\nabla _{f}|_{{\cal Q}({\Bbb P}^{1},A)}$ are injective; thus, 
\begin{eqnarray*}
\dim {\cal D}^{\flat }(f,A) &\leq &\dim {\cal Q}({\Bbb P}^{1},A^{+})/\nabla
_{f}{\cal Q}({\Bbb P}^{1},A)\,=\,\dim {\cal Q}({\Bbb P}^{1},A^{+})-\dim
\nabla _{f}{\cal Q}({\Bbb P}^{1},A) \\
&=&\dim {\cal Q}({\Bbb P}^{1},A^{+})-\dim {\cal Q}({\Bbb P}%
^{1},A)\,=\,\#(A^{+}-A)\,=\,\delta (f)
\end{eqnarray*}
so $\gamma (f)=\sup_{A}\sum_{\langle x\rangle \subseteq A}\gamma _{\langle
x\rangle }=\sup_{A}\dim {\cal D}^{\flat }(f,A)\leq \delta (f)$ by
Proposition \ref{dimd}. $\square $

\section{Local Considerations\label{local}}

Here we prove Proposition \ref{dimd}: that $\dim {\cal D}_{\langle x\rangle
}^{\flat }(f)=\gamma _{\langle x\rangle }$ for each cycle $\langle x\rangle $
of a rational map $f$. Note that ${\cal D}_{\langle x\rangle }(f)$ depends
only the local behavior of $f$ along $\langle x\rangle $: indeed, it follows
from (\ref{ordpush}) that if $q\in {\cal M}({\Bbb P}^{1})$ with $%
[f_{*}q]_{\langle x\rangle }=\lambda [q]_{\langle x\rangle }\neq 0$ then $%
[q]_{w}=0$ for every $w\notin \langle x\rangle $ in the backward orbit of $%
\langle x\rangle $. It similarly follows that ${\rm ord}_{x}q=-2$ if $%
\langle x\rangle $ is superattracting, in which case $\lambda =\frac{1}{\deg
_{x}f^{\kappa }}$ so that $[f_{*}q]_{\langle x\rangle }\neq [q]_{\langle
x\rangle }$; on the other hand, if $\rho \neq 0$ then $[f_{*}q]_{\langle
x\rangle }=\lambda [q]_{\langle x\rangle }$ precisely when $%
[f^{*}q]_{\langle x\rangle }=\frac{1}{\lambda }[q]_{\langle x\rangle }$.

\begin{lemma}
\label{algdiv}The space ${\cal D}_{\langle x\rangle }(f)$ of invariant
divergences is computed in terms of the formal invariants of $\langle
x\rangle $:

\begin{itemize}
\item  If $\langle x\rangle $ is superattracting then ${\cal D}_{\langle
x\rangle }(f)=0$.

\item  If $\langle x\rangle $ is attracting, repelling or irrationally
indifferent then ${\cal D}_{\langle x\rangle }(f)$ is the 
\mbox{$1$-dimensional} space generated by $\frac{d\zeta
^{2}}{\zeta ^{2}}$, for any local coordinate $\zeta $ vanishing at $x$.

\item  If $\langle x\rangle $ is parabolic then ${\cal D}_{\langle x\rangle
}(f)$ is the direct sum of the $\nu $-dimensional subspace ${\cal D}%
_{\langle x\rangle }^{\circ }(f)$ generated by 
\[
\frac{d\zeta ^{2}}{\zeta ^{2}},...\,,\frac{d\zeta ^{2}}{\zeta ^{\ell n+2}}%
,...\,,\frac{d\zeta ^{2}}{\zeta ^{N-n+2}}
\]
and the \mbox{$1$-dimensional} subspace generated by 
\[
\frac{d\zeta ^{2}}{(\zeta ^{N+1}-\beta \zeta ^{2N+1})^{2}}=\frac{d\zeta ^{2}%
}{\zeta ^{2N+2}}+2\beta \frac{d\zeta ^{2}}{\zeta ^{N+2}}+3\beta ^{2}\frac{%
d\zeta ^{2}}{\zeta ^{2}}+O\left( \frac{d\zeta ^{2}}{\zeta }\right) ,
\]
for any local coordinate $\zeta $ as in {\em (\ref{parnf})}.
\end{itemize}
\end{lemma}

\noindent {\bf Proof: }In view of the discussion above, it suffices to
determine when $[f^{*}q]_{\langle x\rangle }=[q]_{\langle x\rangle }$.
Assume without loss of generality that $x$ is a fixed point with $\rho \neq
0 $, and let $\zeta $ be any local coordinate vanishing at $x$; as $f^{*}%
\frac{d\zeta ^{2}}{\zeta ^{j+2}}=\frac{\left( \rho +O(\zeta )\right) ^{2}}{%
\left( \rho \zeta +O(\zeta ^{2})\right) ^{j+2}}\,d\zeta ^{2}=\rho ^{-j}\frac{%
d\zeta ^{2}}{\zeta ^{j+2}}+O\left( \frac{d\zeta ^{2}}{\zeta ^{j+1}}\right) $
for any integer $j$, it follows that ${\cal D}_{x}(f)={\Bbb C[}\frac{d\zeta
^{2}}{\zeta ^{2}}]_{x}$ unless $\rho $ is a root of unity. Suppose now that 
$x$ is parabolic, and let $\zeta $ be a local coordinate as in (\ref{parnf});
then 
\[
\left\{ 
\begin{array}{lll}
f^{*}\zeta ^{k} & = & \left( \rho \left( \zeta +\zeta ^{N+1}+\left( \tfrac{%
N+1}{2}-\beta \right) \zeta ^{2N+1}+O(\zeta ^{2N+2})\right) \right) ^{k} \\ 
& = & \rho ^{k}\zeta ^{k}\left( 1+k\zeta ^{N}+k\left( \tfrac{N+k}{2}-\beta
\right) \zeta ^{2N}+O(\zeta ^{2N+1})\right) , \\ 
f^{*}d\zeta ^{2} & = & \left( \rho \left( 1+(N+1)\zeta ^{N}+(2N+1)\left( 
\tfrac{N+1}{2}-\beta \right) \zeta ^{2N}+O(\zeta ^{2N+1})\right) d\zeta
\right) ^{2} \\ 
& = & \rho ^{2}(1+(2N+2)\zeta ^{N}+(3N^{2}+5N+2-(4N+2)\beta )\zeta
^{2N}+O(\zeta ^{2N+1}))d\zeta ^{2},
\end{array}
\right. 
\]
so 
\[
f^{*}\tfrac{d\zeta ^{2}}{\zeta ^{j+2}}=\bar{\rho}^{j}\left( 1+(2N-j)\zeta
^{N}+(3N^{2}-\tfrac{5}{2}jN+\tfrac{1}{2}j^{2}+(j-4N)\beta )\zeta
^{2N}+O(\zeta ^{2N+1})\right) \tfrac{d\zeta ^{2}}{\zeta ^{j+2}}. 
\]
In particular, 
\begin{eqnarray*}
{\rm ord}_{x}\left( f^{*}\frac{d\zeta ^{2}}{\zeta ^{j+2}} - 
                  \frac{d\zeta ^{2}}{\zeta ^{j+2}}\right) 
&\geq& N-(j+2) \;\;\mbox{if $n|j$, with equality for $j \neq 2N$, and}\\
{\rm ord}_{x}\left(f^{*}\frac{d\zeta ^{2}}{\zeta ^{j+2}} - 
\frac{d\zeta ^{2}}{\zeta ^{j+2}}\right) &=&-(j+2) \;\mbox{otherwise;}
\end{eqnarray*}
as 
\begin{eqnarray*}
f^{*}\frac{d\zeta ^{2}}{\zeta ^{2N+2}}+2\beta f^{*}\frac{d\zeta ^{2}}{\zeta
^{N+2}} &=&(1-2N\beta \zeta ^{2N})\frac{d\zeta ^{2}}{\zeta ^{2N+2}}+2\beta
(1+N\zeta ^{N})\frac{d\zeta ^{2}}{\zeta ^{N+2}}+O\left( \frac{d\zeta ^{2}}{%
\zeta }\right) \\
&=&\frac{d\zeta ^{2}}{\zeta ^{2N+2}}+2\beta \frac{d\zeta ^{2}}{\zeta ^{N+2}}%
+O\left( \frac{d\zeta ^{2}}{\zeta }\right) ,
\end{eqnarray*}
it follows that ${\cal D}_{x}(f)={\Bbb C}\left[ \frac{d\zeta ^{2}}{(\zeta
^{N+1}-\beta \zeta ^{2N+1})^{2}}\right] _{x}\oplus {\cal D}_{x}^{\circ }(f)$%
. $\square $

\medskip

These bases are actually canonical. Indeed, 
$[\frac{d\zeta ^{2}}{\zeta ^{2}}]_{x}$ is independent of the choice of local
coordinate~$\zeta$, and $[\frac{d\zeta ^{2}}{\zeta ^{\ell n+2}}]_{x}$ is
invariant under coordinate 
changes \mbox{$\zeta \mapsto \rho ^{j}\zeta + O(\zeta ^{\ell n+2})$} for
$j\in {\Bbb Z}$, while 
\mbox{$[q_{f}]_{x}=\left[ \frac{d\zeta ^{2}}
                      {(\zeta ^{N+1}-\beta \zeta ^{2N+1})^{2}}\right] _{x}$}
is invariant under the coordinate changes respecting the
normal form (\ref{parnf}); more precisely, consideration of the change of
variable $Z^{N}=\tau \zeta ^{N}$ shows that $[q_{f}]_{x}\in {\cal D}%
_{\langle x\rangle }(F)$ for any germ 
\[
F(\zeta )=e^{2\pi ik/N}\left( \zeta +\tau \zeta ^{N+1}+\left( \tfrac{N+1}{2}%
\tau ^{2}-\beta \right) \zeta ^{2N+1}+O(\zeta ^{2N+2})\right) 
\]
with $(k,\tau )\in {\Bbb Z}\times {\Bbb C}$. As 
\begin{equation}
f^{m}(\zeta )=\rho ^{m}\left( \zeta +m\zeta ^{N+1}+\left( \tfrac{N+1}{2}%
m^{2}-\beta \right) \zeta ^{2N+1}+O(\zeta ^{2N+2})\right)  \label{iterate}
\end{equation}
it similarly follows that $[q_{f}]_{x}=\frac{1}{m^{2}}[q_{f^{m}}]_{x}$ for
any integer $m$. In fact, there is always a unique normalized {\em formal}
quadratic differential $q_{f}$ with $f^{*}q_{f}=q_{f}$. If $x$ is formally
linearizable then $q_{f}=\frac{d\zeta ^{2}}{\zeta ^{2}}$ for any choice of
formal linearizing coordinate $\zeta $. Otherwise, $x$ is parabolic and $%
q_{f}=\frac{1}{n^{2}}q_{f^{n}}=\frac{1}{n^{2}}\eta _{f^{n}}^{2}$ where $\eta
_{f^{n}}$ is the formal linear differential dual to the unique formal vector
field $v_{f^{n}}$ whose formal exponential is~$f^{n}$: necessarily, $%
v_{f^{n}}\nolinebreak =\nolinebreak n\left( \zeta ^{N+1}\nolinebreak
-\nolinebreak \beta \zeta ^{2N+1}\nolinebreak +\nolinebreak O(\zeta
^{2N+2})\right) \frac{\partial }{\partial \zeta }$ in any local coordinate
as in (\ref{parnf}). We pursue this more intrinsic approach in \cite{algdyn}.

Recall that if $\xi $ is a smooth vector field and $\varpi $ is a smooth
2-form, then the flux across an oriented smooth curve $\Gamma $ is $%
\int_{\Gamma }\iota _{\xi }\varpi $, where $\iota _{\xi }$ is the interior
product with respect to $\xi $. In particular, if $\varpi $ has an isolated
singularity at $x$ then there is a flux across $\partial U$ for any
sufficiently small smoothly bounded neighborhood; in view of Cartan's
Formula ${\cal L}_{\xi }\vartheta =\iota _{\xi }d\vartheta +d\iota _{\xi
}\vartheta $, it follows by Stokes Theorem that there is a well-defined
asymptotic flux 
\[
\int_{\partial U}\iota _{\xi }\varpi -\int_{\partial U}{\cal L}_{\xi }\varpi
=\lim_{U\searrow x}\int_{\partial U}\iota _{\xi }\varpi 
\]
from any singularity where the Lie derivative ${\cal L}_{\xi }\varpi $
remains integrable. Similarly, if $[q]_{\langle x\rangle }\in {\cal D}%
_{\langle x\rangle }(f)$ then the local action of $f_{*}$ creates or
destroys a well-defined quantity of mass. Indeed, the net \hbox{$q$-mass}
exiting $U\supset \langle x\rangle $ is $\int_{f(U)-U}|q|-\int_{U-f(U)}|q|$;
moreover, if $\tilde{U}\supset \langle x\rangle $ is contained in $U\cap f(U)
$ then \mbox{$\int_{f(U)-U}|q|-\int_{U-f(U)}|q|=\int_{f(U)-\tilde{U}%
}|q|\,-\,\int_{U-\tilde{U}}|q|$} and
\mbox{$\int_{f(U)-f(\tilde{U})}|q|\,-\,\int_{U-% 
\tilde{U}}|q|=\int_{U-\tilde{U}}(f^{*}|q|-|q|)$}, hence 
\[
\left( \int_{f(U)-U}|q|-\int_{U-f(U)}|q|\right) -\left( \int_{f(\tilde{U})-%
\tilde{U}}|q|\,-\,\int_{\tilde{U}-f(\tilde{U})}|q|\right) =\int_{U-\tilde{U}%
}(f^{*}|q|-|q|),
\]
so that the quantity 
\mbox{$\int_{f(U)-U}|q|\,-\,\int_{U-f(U)}|q|\,-\,\int_{U}(f^{*}|q|-|q|)$} is
independent of $U$. This quantity makes sense because 
\mbox{$\left|\,f^{*}|q|-|q|\,\right| \leq |f^{*}q-q|$} is
locally integrable.  We set
\begin{eqnarray*}
{\rm Res}_{\langle x\rangle }(f:q) &=&\frac{1}{2\pi }\left(
\int_{f(U)-U}|q|\,-\,\int_{U-f(U)}|q|\,-\,\int_{U}(f^{*}|q|-|q|)\right)  \\
&=&\frac{1}{2\pi }\lim_{U\searrow \langle x\rangle }\left(
\int_{f(U)-U}|q|\,-\,\int_{U-f(U)}|q|\right) ;
\end{eqnarray*}
as $\left| \,|q|-|\hat{q}|\,\right| \leq |q-\hat{q}|$ is locally integrable
for quadratic differentials $q$ and $\hat{q}$ with the same algebraic
divergence, the invariant ${\rm Res}_{\langle x\rangle }(f:q)$ depends only
on $[q]_{\langle x\rangle }$. Note that ${\rm Res}_{x}(f^{m}:q)=m{\rm Res}%
_{x}(f:q)$ for a fixed point $x$ and any integer $m$, because 
\[
\int_{f^{m}(U)-U}|q|-\int_{U-f^{m}(U)}|q|=\sum_{j=0}^{m-1}\left(
\int_{f^{j+1}(U)-f^{j}(U)}|q|-\int_{f^{j}(U)-f^{j+1}(U)}|q|\right) ;
\]
similarly, ${\rm Res}_{\langle x\rangle }(f:q)={\rm Res}_{x}(f^{\kappa }:q)$
for any $\kappa $-cycle $\langle x\rangle $, as we may take $U$ to be the
disjoint union $\bigcup_{j=0}^{\kappa -1}f(U_{x})$ for an appropriately
small neighborhood $U_{x}\ni x$.

\begin{lemma}
\label{metdiv}The dynamical residue ${\rm Res}_{\langle x\rangle }(f:q)$ is
computed in terms of the formal invariants of $\langle x\rangle $:

\begin{itemize}
\item  If $\langle x\rangle $ is attracting, indifferent or repelling then $%
{\rm Res}_{\langle x\rangle }(f:q)=|c|\log |\rho |$ for the invariant
divergence $[q]_{\langle x\rangle }=c[\frac{d\zeta ^{2}}{\zeta ^{2}}%
]_{\langle x\rangle }$.

\item  If $\langle x\rangle $ is parabolic then ${\rm Res}_{\langle x\rangle
}(f:q)=|c|\Re \beta $ for any invariant divergence $[q]_{\langle x\rangle }$
in the hyperplane $c[q_{f}]_{\langle x\rangle }+{\cal D}_{\langle x\rangle
}^{\circ }(f)$.
\end{itemize}
\end{lemma}

\noindent {\bf Proof: }Assume without loss of generality that $x$ is a fixed
point, and let $\zeta $ be a local coordinate vanishing at $x$. The choice
of $U\ni x$ is immaterial so we may further assume that $\partial U$ is
smooth; it follows by Stokes Theorem that 
\begin{eqnarray*}
{\rm Res}_{x}(f:q) &=&\frac{1}{2\pi }\left( \int_{\partial f(U)}\vartheta
\,-\,\int_{\partial f(U)}\vartheta \,-\,\int_{U}(f^{*}|q|-|q|)\right) \\
&=&\frac{1}{2\pi }\left( \int_{\partial U}(f^{*}\vartheta -\vartheta
)\,-\,\int_{U}(f^{*}|q|-|q|)\right) \,=\,\frac{1}{2\pi }\lim_{U\searrow
x}\int_{\partial U}(f^{*}\vartheta -\vartheta )
\end{eqnarray*}
for any 1-form $\vartheta $ with $|q|=d\vartheta $. Let $t\mapsto f_{t}$ be
a smooth path of holomorphic germs connecting the identity $f_{0}$ to $%
f=f_{1}$, and consider the holomorphic vector fields $v_{t}$ with $\frac{d}{%
dt}f_{t}=v_{t}\circ f_{t}$. Then 
\begin{eqnarray*}
f^{*}\vartheta -\vartheta &=&\int_{0}^{1}\left( \frac{d}{dt}%
f_{t}^{*}\vartheta \right) dt=\int_{0}^{1}f_{t}^{*}({\cal L}_{\xi
_{t}}\vartheta )\,dt=\int_{0}^{1}{\cal L}_{f_{t}^{*}\xi
_{t}}(f_{t}^{*}\vartheta )\,dt \\
&=&\int_{0}^{1}\iota _{f_{t}^{*}\xi _{t}}(f_{t}^{*}|q|)\,dt\,+\,d\left(
\int_{0}^{1}\iota _{f_{t}^{*}\xi _{t}}(f_{t}^{*}\vartheta )\,dt\right)
\end{eqnarray*}
by Cartan's Formula, where $\xi _{t}=2\Re v_{t}$ corresponds to $v_{t}$
under the standard identification of real and complex tangent spaces; as $%
\iota _{f_{t}^{*}\xi _{t}}(f_{t}^{*}|q|)=\Re \iota _{2f_{t}^{*}v_{t}}(\Re
f_{t}^{*}|q|)$ and $\Re f_{t}^{*}|q|=f_{t}^{*}|q|$, it follows that 
\[
{\rm Res}_{x}(f:q)=\frac{1}{\pi }\lim_{U\searrow x}\Re \int_{\partial
U}\int_{0}^{1}\iota _{f_{t}^{*}v_{t}}(f_{t}^{*}|q|)\,dt. 
\]
In particular, if $f_{t}^{*}v_{t}=v+O\left( \zeta ^{k}\frac{\partial }{%
\partial \zeta }\right) $ for some $k\geq -{\rm ord}_{x}q$ and if $%
f_{t}^{*}q=q+O\left( \frac{d\zeta ^{2}}{\zeta }\right) $, then 
\[
\int_{0}^{1}\iota _{f_{t}^{*}v_{t}}(f_{t}^{*}|q|)\,dt=\int_{0}^{1}\left(
\iota _{v}|q|+O(d\bar{\zeta})\right) dt=\iota _{v}|q|+O(d\bar{\zeta}),\]
so that 
\[
{\rm Res}_{x}(f:q)=\frac{1}{\pi }\lim_{U\searrow x}\Re \int_{\partial
U}\left( \iota _{v}|q|+O(d\bar{\zeta})\right) =\Re \lim_{U\searrow x}\frac{1%
}{\pi }\int_{\partial U}\iota _{v}|q|=\Re \lim_{U\searrow x}\frac{1}{\pi }%
\int_{\partial U}\overline{\iota _{v}|q|}. 
\]

Suppose first that $q=c\frac{d\zeta ^{2}}{\zeta ^{2}}$, and take $%
f_{t}(\zeta )=e^{t\log \rho }\zeta +t\left( f(\zeta )-\rho \zeta \right) $
for some choice of $\log \rho $. Then $f_{t}(\zeta )=e^{t\log \rho }\zeta
+O(\zeta ^{2})$, so 
\[
f_{t}^{*}v_{t}=\left. \frac{df_{t}}{dt}\right/ \frac{df_{t}}{d\zeta }=\left( 
\frac{(\log \rho )e^{t\log \rho }\zeta +O(\zeta ^{2})}{e^{t\log \rho
}+O(\zeta )}\right) \frac{\partial }{\partial \zeta }=v+O\left( \zeta ^{2}%
\frac{\partial }{\partial \zeta }\right) 
\]
where $v=(\log \rho )\zeta \frac{\partial }{\partial \zeta }$, and $%
f_{t}^{*}q=q+O\left( \frac{d\zeta ^{2}}{\zeta }\right) $ because ${\rm ord}%
_{x}q\geq -2$; as $|d\zeta ^{2}|\nolinebreak =\nolinebreak \frac{i}{2}d\zeta
\nolinebreak \wedge \nolinebreak d\bar{\zeta}$ and therefore $\iota
_{v}|q|=(\log \rho )\zeta \left( \frac{i|c|}{2|\zeta |^{2}}\right) \iota _{%
\frac{\partial }{\partial \zeta }}\left( d\zeta \wedge d\bar{\zeta}\right) =%
\frac{i}{2}|c|(\log \rho )\frac{d\bar{\zeta}}{\bar{\zeta}}$, it follows that 
\[
{\rm Res}_{x}(f:q)=\Re \frac{1}{2\pi i}\oint |c|\overline{(\log \rho )}\frac{%
d\zeta }{\zeta }=|c|\log |\rho |.
\]
Assume now that $x$ is parabolic; let $\zeta $ be a local coordinate as in (%
\ref{parnf}), and take 
\[
f_{t}(\zeta )=\zeta +t\left( f^{n}(\zeta )-\zeta +%
\tfrac{N+1}{2}(t-1)n^{2}\zeta ^{2N+1}\right)
\]
so that $f_{1}=f^{n}$.
Then $f_{t}(\zeta )=\zeta +tn\zeta ^{N+1}+tn\left( \tfrac{N+1}{2}tn-%
\beta \right)\zeta ^{2N+1}+O(\zeta ^{2N+2})$ by (\ref{iterate}), so 
\begin{eqnarray*}
f_{t}^{*}v_{t} &=&\left. \frac{df_{t}}{dt}\right/ \frac{df_{t}}{d\zeta }%
=\left( \frac{n\zeta ^{N+1}+n\left( (N+1)tn-\beta \right) \zeta
^{2N+1}+O(\zeta ^{2N+2})}{1+tn(N+1)\zeta ^{N}+O(\zeta ^{2N})}\right) \frac{%
\partial }{\partial \zeta } \\
&=&v+O\left( \zeta ^{2N+2}\frac{\partial }{\partial \zeta }\right) 
\end{eqnarray*}
where $v=n(\zeta ^{N+1}-\beta \zeta ^{2N+1})\frac{\partial }{\partial \zeta }
$. By the remarks following Lemma \ref{algdiv}, if 
\mbox{$q= c\frac{d\zeta ^{2}}{(\zeta ^{N+1}-\beta \zeta ^{2N+1})^{2}}%
+O\left( \frac{d\zeta ^{2}}{\zeta ^{N+1}}\right) $} then 
\mbox{$f_{t}^{*}q=q+O\left( \frac{d\zeta ^{2}}{\zeta }\right) $}; as 
\begin{eqnarray*}
\iota _{v}|q| &=&n\left( \zeta ^{N+1}-\beta \zeta ^{2N+1}\right) \left( 
\frac{i}{2}|c|\frac{1+O(\zeta ^{N+1})}{|\zeta ^{N+1}-\beta \zeta ^{2N+1}|^{2}%
}\right) \iota _{\frac{\partial }{\partial \zeta }}\left( d\zeta \wedge d%
\bar{\zeta}\right)  \\
&=&\frac{i}{2}n|c|\left( \frac{1}{\overline{\zeta ^{N+1}-\beta \zeta ^{2N+1}}%
}+O(1)\right) d\bar{\zeta}=\frac{i}{2}n|c|\left( \frac{d\bar{\zeta}}{\bar{%
\zeta}^{N+1}}+\bar{\beta}\frac{d\bar{\zeta}}{\bar{\zeta}}\right) +O(d\bar{%
\zeta}),
\end{eqnarray*}
it follows that ${\rm Res}_{x}(f:q)=\frac{1}{n}{\rm Res}_{x}(f^{n}:q)=|c|\Re 
\frac{1}{2\pi i}\oint \left( \frac{d\zeta }{\zeta ^{N+1}}+\beta \frac{d\zeta 
}{\zeta }\right) =|c|\Re \beta .$ $\square $

\medskip

Recall that ${\cal D}_{\langle x\rangle }^{\flat }(f)=\{[q]_{\langle
x\rangle }\in {\cal D}_{\langle x\rangle }(f):$ ${\rm Res}_{\langle x\rangle
}(f:q)\leq 0\}$ by definition. The above discussion shows that 
\[
{\cal D}_{\langle x\rangle }^{\flat }(f)=\left\{ 
\begin{array}{ll}
0 & \text{if }\langle x\rangle \text{ is repelling or superattracting} \\ 
{\Bbb C[}\frac{d\zeta ^{2}}{\zeta ^{2}}]_{\langle x\rangle } & \text{if }%
\langle x\rangle \text{ is attracting or irrationally indifferent} \\ 
{\cal D}_{\langle x\rangle }^{\circ }(f) & \text{if }\langle x\rangle \text{
is parabolic-repelling} \\ 
{\cal D}_{\langle x\rangle }(f) & \text{if }\langle x\rangle \text{ is
parabolic-attracting or parabolic-indifferent,}
\end{array}
\right. 
\]
and Proposition \ref{dimd} now follows from the definition of $\gamma
_{\langle x\rangle }$.

\section{Global Considerations \label{global}}

We now prove Proposition \ref{main}: the injectivity of $\nabla _{f}:{\cal Q}%
^{\flat }(f)\rightarrow {\cal Q}({\Bbb P}^{1})$. Recall that $||f_{*}q||\leq
||q||$ for any quadratic differential $q$, but note that this inequality is
vacuous when $||q||=\infty $; nevertheless, we may still identify the mass
decrease due to cancellation as 
\[
{\rm Dec}(f:q)\nolinebreak =\nolinebreak \int_{{\Bbb P}%
^{1}}(f_{*}|q|-|f_{*}q|) 
\]
which is always nonnegative and might be infinite. Clearly, ${\rm Dec}%
(f:q)\nolinebreak =\nolinebreak ||q||\nolinebreak -\nolinebreak ||f_{*}q||$
for any integrable $q$, so the following extension of Thurston's Contraction
Principle already shows injectivity on ${\cal Q}({\Bbb P}^{1})$:

\begin{lemma}
\label{contraction}Let $f:{\Bbb P}^{1}\rightarrow {\Bbb P}^{1}$ be a
rational map of degree $D>1$; assume that $f$ is not a Latt\`{e}s example,
and let $q\neq 0$ be a meromorphic quadratic differential in $\ker \nabla
_{f}$. Then ${\rm Dec}(f:q)>0$.
\end{lemma}

\noindent {\bf Proof: }Fix an open disk $U\subset {\Bbb P}^{1}-S(f)$ and an
inverse branch $\phi :U\rightarrow {\Bbb P}^{1}$. If ${\rm Dec}(f:q)=0$ then 
$|\phi ^{*}q\nolinebreak +\nolinebreak \psi ^{*}q|\nolinebreak =\nolinebreak
|\phi ^{*}q|\nolinebreak +\nolinebreak |\psi ^{*}q|$ for any inverse branch $%
\psi \nolinebreak :\nolinebreak U\nolinebreak \rightarrow \nolinebreak {\Bbb %
P}^{1}$, so the meromorphic function $\lambda _{\psi }\nolinebreak =\frac{%
\psi ^{*}q}{\phi ^{*}q}$ is almost everywhere real and positive. As a
real-valued meromorphic function is necessarily constant, $f_{*}q=\lambda
\phi ^{*}q$ on $U$, so $f^{*}f_{*}q\nolinebreak =\nolinebreak \lambda
f^{*}\phi ^{*}q\nolinebreak =\nolinebreak \lambda q$ on $\phi (U)$ and
therefore globally, where $\lambda =\sum_{\psi }\lambda _{\psi }\in {\Bbb R}%
^{+}$; it follows from (\ref{pushpull}) that $Df_{*}q=f_{*}f^{*}f_{*}q=%
\lambda f_{*}q$, so in fact $\lambda =D$.

As $\nabla _{f}q=0$, it further follows that $f^{*}q=f^{*}f_{*}q=Dq$, and we
claim that this is impossible unless $f$ is a Latt\`{e}s example. Indeed, if 
$f^{*}q=\lambda q$ for any $\lambda \neq 0$ then {\rm ord}$_{x}f^{*}q=\,$%
{\rm ord}$_{x}q$ for each point $x$ on the Riemann sphere. In view of (\ref
{qorder}), the finite set 
\[
\Phi  = \{x \in {\Bbb P}^{1}: {\rm ord}_{x}q \leq -2 \;\mbox{or}\; {\rm
ord}_{x}q \geq 1\} 
\]
is backward invariant, hence $\#\Phi \leq 2$ with any $x\in \Phi $
superattracting of period 1 or 2 (see \cite{BH} or \cite{Mil}); but then 
{\rm ord}$_{x}q=-2$ which is only possible for $\lambda =D^{2}$, so $\Phi
=\emptyset $. In particular, $q$ is nowhere vanishing with only simple
poles; in fact, there are precisely four such poles, as $\sum_{x\in {\Bbb P}%
^{1}}${\rm ord}$_{x}q=-4$ for any quadratic differential on ${\Bbb P}^{1}$.
It now follows from (\ref{qorder}) that the set of poles is forward
invariant, that every preimage of a pole is either a pole or a simple
critical point, that every critical value is a pole, and that no critical
point is a pole: this is precisely the description of the Latt\`{e}s
examples given in \cite{DH2}. $\square $

\medskip

To complete the proof of Proposition \ref{main}, we show $\,{\rm Dec}(f:q)=0$
for any quadratic differential $q$ in ${\cal Q}^{\flat }(f)\nolinebreak \cap
\nolinebreak \ker \nabla _{f}$. Consider the total dynamical residue 
\[
{\rm Res}(f\nolinebreak :\nolinebreak q)\nolinebreak =\nolinebreak
\sum_{\langle x\rangle \subset {\Bbb P}^{1}}{\rm Res}_{\langle x\rangle
}(f\nolinebreak :\nolinebreak q);
\]
note that this quantity is defined whenever $q\in {\cal Q}(f)$, and that $%
{\rm Res}(f:q)\leq 0$ for any $q\in {\cal Q}^{\flat }(f)$ As ${\rm Dec}%
(f:q)\geq 0$ for any $q\in {\cal M}({\Bbb P}^{1})$, the desired conclusion
follows from the Balance Principle: that 
\[
{\rm Dec}(f:q)=2\pi {\rm Res}(f:q)
\]
for every $q\in \ker \nabla _{f}$. This identity is an immediate consequence
of the following:

\begin{lemma}
\label{lowerbd}Let $f:{\Bbb P}^{1}\rightarrow {\Bbb P}^{1}$ be a rational
map, and let $q$ be a meromorphic quadratic differential in ${\cal Q}(f)$.
Then $||\nabla _{f}q||\geq |{\rm Dec}(f\nolinebreak :\nolinebreak
q)\nolinebreak -2\pi \,{\rm Res}(f\nolinebreak :\nolinebreak q)|$.
\end{lemma}

\noindent {\bf Proof:} Recall that $E=\{x\in {\Bbb P}^{1}:{\rm ord}_{x}q\leq
-2\}$ consists of finitely many cycles, none superattracting. Let $U_{x}\ni x
$ be pairwise disjoint open disks in ${\Bbb P}^{1}-S(f)$, and set $%
U\nolinebreak =\nolinebreak \bigcup_{x\in E}U_{x}$; we may arrange that $%
U\cap f^{-1}(U)=\bigcup_{x\in E}\left( U_{x}\cap \phi (U_{f(x)})\right) $
where $\phi :U\rightarrow {\Bbb P}^{1}$ with $\phi (E)=E$ consists of the
distinguished inverse branches of $f$.

We claim that $f_{*}|q|-|f_{*}q|$ is integrable, so that ${\rm Dec}%
(f\nolinebreak :\nolinebreak q)<\infty $: indeed, both $f_{*}|q|$ and $%
|f_{*}q|$ are integrable on ${\Bbb P}^{1}-U$, and the restriction of $%
f_{*}|q|-|f_{*}q|$ to $U$ is 
\[
(\phi ^{*}|q|-|f_{*}q|)+\sum_{\psi \neq \phi }\psi ^{*}|q|
\]
where $\left| \,\phi ^{*}|q|-|f_{*}q|\,\right| \leq \left| \sum_{\psi \neq
\phi }\psi ^{*}q\right| \leq \sum_{\psi \neq \phi }\psi ^{*}|q|$. Moreover, $%
\left| \,|q|\nolinebreak -\nolinebreak |f_{*}q|\,\right| \nolinebreak \leq
\nolinebreak |q\nolinebreak -\nolinebreak f_{*}q|$ so $|q|-|f_{*}q|$ is
integrable, and it follows that $|q|-f_{*}|q|$ is also integrable: in fact, 
\begin{eqnarray*}
\int_{{\Bbb P}^{1}-U}\left( |q|-f_{*}|q|\right) \, &=&\,\int_{{\Bbb P}%
^{1}-U}|q|\,-\,\int_{{\Bbb P}^{1}-f^{-1}(U)}|q|\, \\
\, &=&\,\int_{f^{-1}(U)-\phi (U)}|q|\,+\,\int_{\phi
(U)-U}|q|\,-\,\int_{U-\phi (U)}|q|
\end{eqnarray*}
and 
\begin{eqnarray*}
\int_{U}\left( |q|-f_{*}|q|\right) \, &=&\,\int_{U}\left( |q|-\phi
^{*}|q|\right) \,-\,\int_{U}\left( f_{*}|q|-\phi ^{*}|q|\right)  \\
\, &=&\,\int_{U}\left( |q|-\phi ^{*}|q|\right) \,-\,\int_{f^{-1}(U)-\phi
(U)}|q|
\end{eqnarray*}
so 
\begin{eqnarray*}
\int_{{\Bbb P}^{1}}\left( |q|-f_{*}|q|\right) \, &=&\,\int_{\phi
(U)-U}|q|\,-\,\int_{U-\phi (U)}|q|\,-\,\int_{U}\left( \phi
^{*}|q|-|q|\right)  \\
\, &=&\,2\pi \sum_{\langle x\rangle \subseteq E}{\rm Res}_{\langle x\rangle
}(\phi :q)\,=\,-2\pi \,{\rm Res}(f:q)\text{.}
\end{eqnarray*}
Consequently, 
\[
||\nabla _{f}q||\,\geq \,\left| \int_{{\Bbb P}^{1}}\left(
f_{*}|q|-|f_{*}q|\right) \,+\,\int_{{\Bbb P}^{1}}\left( |q|-f_{*}|q|\right)
\right| \,=\,\left| {\rm Dec}(f:q)-2\pi \,{\rm Res}(f:q)\right| \text{.}
\]
$\square $

\bigskip \noindent {\bf Address:}

\noindent Department of Mathematics

\noindent Cornell University

\noindent Ithaca, NY 14853-7901

\medskip \noindent {\bf E-mail: } adame@math.cornell.edu

\end{document}

%% file: imsmark.tex
\def\IMSmarkvadjust{0 pt}
\def\IMSmarkhadjust{0 pt}
\def\SBIMSMark#1#2#3{
 \font\SBF=cmss10 at 10 true pt
 \font\SBI=cmssi10 at 10 true pt
 \setbox0=\hbox{\SBF Stony Brook IMS Preprint \##1}
 \setbox2=\hbox to \wd0{\hfil \SBI #2}
 \setbox4=\hbox to \wd0{\hfil \SBI #3}
 \setbox6=\hbox to \wd0{\hss
             \vbox{\hsize=\wd0 \parskip=0pt \baselineskip=10 true pt
                   \copy0 \break%
                   \copy2 \break% 
                   \copy4 \break}}
 \dimen0=\ht6   \advance\dimen0 by \vsize \advance\dimen0 by 8 true pt
                \advance\dimen0 by -\pagetotal
	        \advance\dimen0 by \IMSmarkvadjust
 \dimen2=\hsize \advance\dimen2 by .25 true in
	        \advance\dimen2 by \IMSmarkhadjust

%
%   Check for publication info
%
%  \newread\jref
  \openin2=publishd.tex
  \ifeof2\setbox0=\hbox to 0pt{}
  \else 
     \setbox0=\hbox to 3.1 true in{
                \vbox to \ht6{\hsize=3 true in \parskip=0pt  \noindent  
                {\SBI Published in modified form:}\hfil\break
                \input publishd.tex 
                \vfill}}
  \fi
  \closein2
  \ht0=0pt \dp0=0pt
 \ht6=0pt \dp6=0pt
 \setbox8=\vbox to \dimen0{\vfill \hbox to \dimen2{\copy0 \hss \copy6}}
 \ht8=0pt \dp8=0pt \wd8=0pt
 \copy8
 \message{*** Stony Brook IMS Preprint #1, #2. #3 ***}
}